\documentclass[11pt]{amsart}

\usepackage{amssymb}
\usepackage{graphicx}
\usepackage{asymptote}

\newcommand{\p}{\partial}

\newcommand{\bo}{\partial M}

\newcommand{\ti}{\tilde}
\newcommand{\ra}{\rightarrow}
\renewcommand{\phi}{\varphi}

\newtheorem{theorem}{Theorem}
\newtheorem{proposition}{Proposition}
\newtheorem{lemma}{Lemma}
\newtheorem{definition}{Definition}
\newtheorem{corollary}{Corollary}
\newtheorem{conjecture}{Conjecture}

\title{A proof of Lens Rigidity in the category of Analytic Metrics}

\author[J. Vargo]{James Vargo}
\address{Department of Mathematics, University of Washington, Seattle, WA 98195}

\begin{document}

\begin{abstract}
Consider a compact Riemannian manifold with boundary.  If all
maximally extended geodesics intersect the boundary at both ends,
then to each geodesic $\gamma(t)$ we can form the triple
$(\dot\gamma(0),\dot\gamma(T),T)$, consisting of the initial and
final vectors of the segment as well as the length between them. The
collection of all such triples comprises the lens data. In this
paper, it is shown that in the category of analytic Riemannian
manifolds, the lens data uniquely determine the metric up to
isometry. There are no convexity assumptions on the boundary, and
conjugate points are allowed, but with some restriction.
\end{abstract}

\maketitle

\section{An introduction including the result proved}
Let $(M,g)$ be a compact, Riemannian manifold with boundary $\bo$,
and let it be non-trapping. That means all geodesics, when maximally
extended, terminate at the boundary at both their ends. Let $SM$
denote its sphere bundle. Then for any vector $v\in \p SM$, the
geodesic $\gamma_v$ originating at $v$ eventually leaves the
manifold after some distance $T$.  Let $\ell(v)$ denote the length
of the geodesic, and let $\Sigma(v)=\dot \gamma_v(T)$ denote its
terminal vector. \[\Sigma:\p SM\rightarrow \p SM\] is called the
scattering map. Together, $\Sigma$ and $\ell$ comprise the lens data
of $(M,g)$.

\begin{figure}[htb]
\centering
\includegraphics{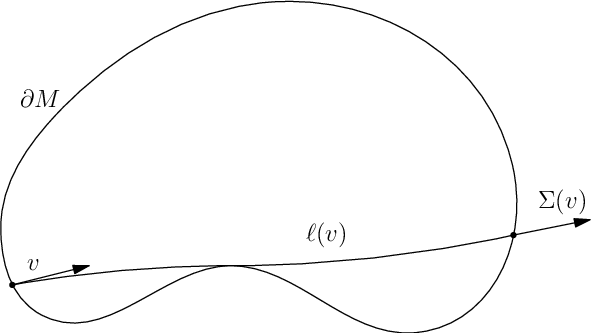}
\caption{The Lens Data}
\end{figure}

The lens rigidity conjecture states that one may recover a
Riemannian manifold up to isometry from its lens data
$(\Sigma,\ell)$. To be more precise, suppose we have two
non-trapping Riemannian manifolds $(M_i,g_i), i=1,2$ which share the
same boundary.  That is, $\bo_1=\bo_2$ (henceforth both will be
denoted $\bo$).  Then for any point $x\in\bo$, there is a natural
correspondence \[\Lambda_x:\p S_xM_1\ra\p S_xM_2.\] Indeed, a unit
vector at the boundary of a Riemannian manifold is uniquely
characterized by its inward normal component and the direction of
its tangential projection. So we shall say that $v_2=\Lambda(v_1)$
if these two quantities agree for $v_1$ and $v_2$, respectively.

\begin{definition}Let $(\Sigma_i,\ell_i)$ denote the lens data for
the manifold $(M_i,g_i): i=1,2$. We shall say that the two manifolds
have the same lens data if
$\Lambda\circ\Sigma_1=\Sigma_2\circ\Lambda$ and $\ell_1=\ell_2\circ
\Lambda$.\end{definition}

\begin{conjecture} \label{cnj lrc}
If $(M_1,g_1)$ and $(M_2,g_2)$ are non-trapping with the same lens
data, then the two manifolds are related by an isometry that fixes
the points of the boundary $\bo$.  That is, there exists a
diffeomorphism $\phi:M_1\rightarrow M_2$ satisfying $\phi|_{\bo}=id$
and $\phi^*g_2=g_1$.
\end{conjecture}

The lens rigidity problem is a generalization of the boundary
rigidity problem.  In that problem, the initial data is taken to be
the boundary distance function
\[\rho_g:\bo\times\bo\rightarrow \mathbb{R}.\]
$\rho_g(x,y)$ is equal to the length of the shortest curve joining
$x$ to $y$.  Of course, metrics related by an isometry fixing the
boundary will also yield the same boundary distance function.  The
boundary rigidity problem is whether this is the only obstruction to
unique recovery of $g$ from $\rho_g$.

If a metric $g$ has the property that the only other metrics with
the same boundary distance function are isometric to $g$, then $g$
is called boundary rigid.  There are many examples of metrics that
are not boundary rigid.  Indeed, $\rho_g$ only records the lengths
of the shortest paths, and it is not hard to construct metrics for
which the shortest paths do not enter certain open subsets of the
manifold.  To circumvent this problem, the assumption of simplicity
is usually made on the metric.

\begin{definition} The Riemannian manifold $(M,g)$ is
simple, if $\partial M$ is strictly convex with respect to $g$, and
for any $x\in  M$,  the exponential map $\exp_x :\exp^{-1}_x(M) \to
M$ is a diffeomorphism.
\end{definition}

A simple manifold has the property of being geodesically convex.
 That is, every pair of points is connected by a unique geodesic and
 that geodesic is length minimizing.  Topologically, a simple
manifold is a ball.  Michel \cite{M} was the first to conjecture
that simple Riemannian manifolds are boundary rigid. This has been
proved recently in two dimensions \cite{PU}. It has also been proved
for subdomains of Euclidean space \cite{Gr}, for metrics close to
the Euclidean \cite{BI}, and symmetric spaces of negative curvature
\cite{BCG}. In \cite{SU2}, Stefanov and Uhlmann proved a local
boundary rigidity result.  If $g$ belongs to a certain generic set
which includes real-analytic metrics, and $g'$ is sufficiently close
to $g$, then $\rho_g=\rho_{g'}$ implies that $g$ and $g'$ be
isometric. For other local results see \cite{CDS}, \cite{E},
\cite{LSU}, \cite{SU1}. It is shown in \cite{SU3} that the lens
rigidity problem is equivalent to the boundary rigidity problem if
the manifold is simple.

If the manifold is not simple, the lens data carries more
information than the boundary distance function.  Indeed, it
includes the lengths of all geodesics, so in the case that $g$ be
non-trapping, these geodesics pass through every point of the
manifold in every direction. However, if the manifold is trapping,
there are examples in which the lens data is not sufficient to
determine the metric, (see \cite{CK}).  There are not many results
on the lens rigidity problem, but the following are notable.  If a
manifold is lens rigid a finite quotient of it is also lens rigid
\cite{C2}. In \cite{SU3}, Stefanov and Uhlmann generalized their
local result for simple metrics to obtain a local lens rigidity
result.  There are some assumptions on conjugate points and a
topological assumption. Assuming these, if $g$ belongs to a certain
generic set which includes real-analytic metrics, and $g'$ is
another metric with the same lens data that, a-priori, is known to
be sufficiently close to $g$, then $g'$ is isometric to $g$.

In this paper, the following statement is proved.

\begin{theorem} \label{thm 1}
Let $(M_i,g_i), i=1,2$ be non-trapping analytic Riemannian manifolds
with a common, analytic boundary $\p M$.
 Further, assume that in each connected component of $S(\bo_1)$,
there exists $(x_0,\xi_0)$ such that $x_0$ is not conjugate to any
points of $\bo$ that lie along the geodesic $\gamma_{x_0,\xi_0}$.
Then if the two manifolds have the same lens data, there must exist
an analytic diffeomorphism $\phi:M_1\rightarrow M_2$ with
$\phi|_{\bo}=id$ and $\phi^*g_2=g_1$.
\end{theorem}

Note the slightly asymmetric nature of the second hypothesis. This
property is used in the proof of Theorem 2 to guarantee the
possibility of a certain construction on the lens data.  Since $g_2$
has the same lens data, the same construction will work
automatically, even though, a priori, there is no reason why the
condition of the theorem should also be true for $g_2$.

\section{Constructing an isometry on a band about the boundary}

Let $(M,g)$ be a general compact Riemannian manifold with boundary.
Let $\nu$ be the field of inward unit normal vectors at the boundary
$\bo$, and let $x_0 \in \bo$ be a boundary point. Then there is a
small neighborhood $N\subset \bo$ of $x_0$ and a number $\epsilon >
0$ such that the mapping \[\exp_{\nu}:N\times [0,\epsilon)\ra M\]
given by $(x',x^n)\mapsto \exp_{x'}(x^n\nu)$ gives a local
coordinate system. These are the boundary normal coordinates.
Through them, the metric has the form
\[ds^2=g_{\alpha\beta}dx'^{\alpha}dx'^{\beta}+(dx^n)^2,\]
where $\alpha,\beta$ are indices running over the tangential
directions.  Now let $\bar{M}$ be an open manifold slightly
extending $M$ and extend $g$ smoothly to $\bar{M}$ (extend by
analytic continuation in the case that $(M,g)$ is analytic).  By
choosing a smaller $\epsilon$ if necessary, we may now use our
boundary normal coordinates as a coordinate system in $\bar{M}$ by
allowing the coordinate $x^n$ to vary over the set
$(-\epsilon,\epsilon)$.

By compactness, we may choose $\epsilon$ uniformly over the whole
boundary.  We may also select $\epsilon$ sufficiently small so that
our boundary normal coordinates give a global diffeomorphism
\[\exp_{\nu}:\p M \times (-\epsilon,\epsilon)\rightarrow V,\]
where $V$ is a neighborhood of $\p M$ in $\bar{M}$.

To show this, it is only necessary to prove the above mapping
injective.  Around each point of $\p M$, let $N$ be a connected open
neighborhood such that $\exp_{\nu}$ defined on $N\times
(-\epsilon,\epsilon)$ is injective.  By compactness, $\p M$ is
covered by a finite number of such neighborhoods $N_1,\dots, N_m$.
There exists a number $\delta>0$ such that for any two points
$x,y\in \bo$, $x,y$ must belong to a common neighborhood $N_i$ if
$d(x,y)<\delta$.  Take $\epsilon$ to be less than half of $\delta$.

Suppose there are points $x,y\in \bo$ and numbers $s,t\in
(-\epsilon,\epsilon)$ such that
\[\exp_{\nu}(x,s)=\exp_{\nu}(y,t).\] Then by the triangle, inequality, $d(x,y)<s+t<\delta$ which
shows that $x,y$ belong to a common neighborhood $N_i$.  But on
$N_i\times (-\epsilon,\epsilon)$, $\exp_{\nu}$ is injective. Hence
$x=y,s=t$.

 We define the manifold $\ti M$ to be $M\cup U$, where $U$ is a collar defined by:

\[U=\{x:-\epsilon \leq x^n \leq 0\}.\]

\begin{figure}[htb]
\centering
\includegraphics{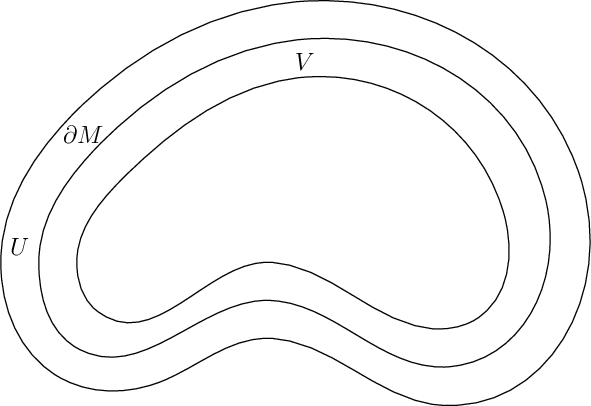}
\caption{$\ti{M}$}

\label{UV}
\end{figure}

Next, note that the set $V$ is a subset of $\ti M$, and can be
described as the set of points in $\ti M$ whose distance from $\bo$
is less than $\epsilon$:

\[V=\{x\in\ti{M}:d(x,\p M)<\epsilon\}.\]
See Figure \ref{UV}.

Theorem 1 relies principally on the following theorem proved by
Stefanov and Uhlmann in \cite{SU3}.

\begin{theorem}  \label{thm_jet}
Let $(M,g)$ be a compact Riemannian manifold with boundary. Let
$(x_0,\xi_0)\in S(\partial M)$ be such that the maximal geodesic
$\gamma_{x_0,\xi_0}$ through it is of finite length,  and assume
that $x_0$ is not conjugate  to any point in $\gamma_{x_0,\xi_0}\cap
\bo$.  Then the jet of $g$ at $x_0$ in boundary normal coordinates
is uniquely determined by the lens data $(\Sigma,\ell)$.
\end{theorem}

\begin{corollary} \label{cor bnc}
Assume $(M,g)$ is analytic with analytic boundary and that, in each
connected component of $S(\bo)$, there is at least one vector
$(x_0,\xi_0)$ satisfying the conditions of the theorem. Then the
lens data uniquely determine the metric $g$ in boundary normal
coordinates.
\end{corollary}

\begin{proof}
As above we let $V$ denote the set of points
$\{x\in\ti{M}:d(x,\bo)\leq\epsilon\}$.  Then by hypothesis, in each
component of $V$ there is at least one point at which the jet of the
metric is determined. Since the metric is analytic, it must be
uniquely determined on all of $V$.
\end{proof}

Now we apply this to our two Riemannian manifolds $(M_i,g_i)$,
taking an $\epsilon$ sufficiently small to work for both. We obtain
$\exp_{\nu_i}:\bo\times[-\epsilon,\epsilon]\ra V_i\subset \ti{M_i}$.
Using these coordinates, both metrics have the form $g_{\alpha
\beta}dx'^{\alpha}dx'^{\beta}+(dx^n)^2$.  By the corollary, the
functions $g_{\alpha\beta}$ coincide for the two metrics throughout
the domains $\bo\times[-\epsilon,\epsilon]$, which means that the
mapping $\phi_0:V_1\ra V_2$ defined by
$\phi_0=\exp_{\nu_2}\circ\exp_{\nu_1}^{-1}$ is an isometry.

Note that $\phi_0|_{\bo}=id$, and $\phi_{0*}(\nu_1)=\nu_2$. In
particular, $\phi_{0*}$ must preserve directions in $T\bo$ and must
preserve components in the normal direction. Thus
$\phi_{0*}|_{SM_1}=\Lambda$, the mapping that relates the lens data
of our two manifolds.\bigskip

\section{Extension of the isometry to the entire manifold}

\begin{figure}[htp]
\centering
\includegraphics{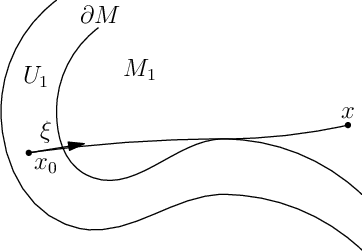}
\hspace{1.5cm}
\includegraphics{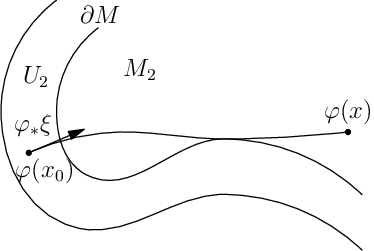}
\caption{} \label{phi}
\end{figure}

The rest of this paper shall be concerned with extending $\phi_0$ to
an isometry $\phi:\ti{M}_1\ra \ti{M_2}$.  If the extension exists,
then it must be uniquely defined. Indeed, given a point $x_0\in U_1$
and a unit vector $\xi$ at $x_0$, we must require

\[
\phi(\exp^{g_1}_{x_0}(t\xi))=\exp^{g_2}_{\phi_0(x_0)}(t\phi_{0*}\xi).
\]
See Figure \ref{phi}.

By the non-trapping assumption, all points $x\in M_1$ lie on a
geodesic originating in $U_1$.  Therefore this equation uniquely
determines the extended mapping $\phi$.  However, it is not at all
clear that the equation yields a well-defined mapping. To get around
this problem, we shall first define a mapping $\ti{\phi}:SM_1\ra
M_2$ and then show that the values of $\ti{\phi}$ only depend on the
basepoint $x\in M_1$.

Choose $(x,v)\in SM_1$, and consider the geodesic $\gamma_{x,-v}$
(note the reversal of $v$).  Let $T_0=T_0(v)\geq 0$ be the time at
which this curve first leaves $M_1$ and enters $U_1$. That is,
\[T_0=T_0(v)=\inf \{t\geq 0:\gamma_{x,-v}(t)\notin M_1\}.\]
This value exists because of the nontrapping assumption. Similarly,
we define $T_1=T_1(v)$ to be the first time after $T_0(v)$ at which
the curve leaves the interior of $U_1$:
\[T_1=T_1(v)=\inf \{t> T_0(v):\gamma_{x,-v}(t)\in \p U_1\}.\]
If, somehow, the curve never leaves the interior of $U_1$, then we
set $T_1=\infty$.

We claim that $T_1>T_0$.  This follows from the assumption that $\p
M$ and the metric are both analytic.  Consequently, a geodesic
cannot have contact of infinite order with the boundary without
being trapped in the boundary. Therefore, at $t=T_0$, we conclude
that there exists $m\geq 1$ for which
\[\begin{array}{cc} \p_t^k (x^n\circ \gamma_{x,-v})(T_0)=0, & 0\leq k <
m;
\\ \p_t^m (x^n\circ \gamma_{x,-v})(T_0)<0. &  \end{array} \]
This shows that for some small amount of time after $T_0$, the
geodesic must remain entirely outside $M_1$.  Hence $T_1>T_0$.

\begin{figure}[htp]
\centering
\includegraphics{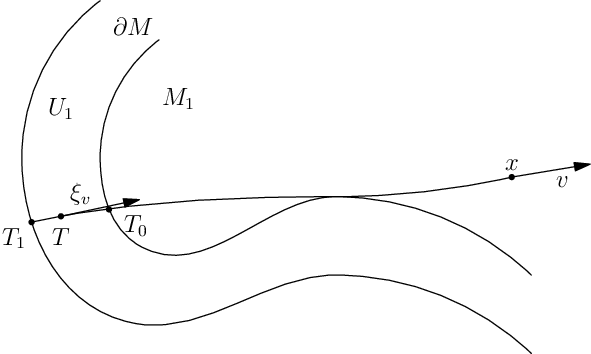}
\caption{The Construction of $\ti{\phi}$} \label{prim}
\end{figure}

 Now we let $T=T(v)$ be an arbitrarily chosen number strictly between
$T_0$ and $T_1$, and we let $\xi_v=-\dot\gamma_{x,-v}(T)$. By
construction, $x=\exp^{g_1}(T\xi_v)$. We define $\ti{\phi}(x,v)$ by:

\[
\ti{\phi}(x,v)=\exp^{g_2}(T\phi_{0*}\xi_v)
\]
See Figure \ref{prim}.

\begin{proposition}
$\ti{\phi}(x,v)$ is a well defined function on $SM_1$ with values in
$M_2$.
\end{proposition}
\begin{proof}
 We must show two things: first, that $\exp^{g_2}(T\phi_{0*}\xi_v)$ is a point in
$M_2$; second, that the value of $\ti{\phi}$ is independent of the
choice of $T$.

The curve $\gamma_{\xi_v}(t)=\exp^{g_1}(t\xi_v):0\leq t\leq T$ is
composed of two segments; the first of which lies within $U_1$, the
second of which lies within $M_1$. The break between the two occurs
at $t=T-T_0$.  We conclude that $\ell(\dot\gamma_{\xi_v}(T-T_0))\geq
T_0$.

The curve $\exp^{g_2}(t\phi_{0*}\xi_v):t\in[0,T]$ is also composed
of two segments of length $T-T_0$ and $T_0$ lying in $U_2$ and $M_2$
respectively.  Indeed, for $t\in [0,T-T_0]$, we have
\begin{equation}\label{corndog}
\exp^{g_2}(t\phi_{0*}\xi_v)=\phi_0(\exp^{g_1}(t\xi_v))\end{equation}
from the fact that $\phi_0$ is an isometry on $U_1$.  Hence the left
side belongs to $U_2$.

To show that the remainder of the curve
$\exp^{g_2}(t\phi_{0*}\xi_v)$ lies in $M_2$, we look at the lens
data.  If we differentiate at $t=T-T_0$, we obtain from equation
\eqref{corndog}, $\phi_{0*}(\dot\gamma_{\xi_v}(T-T_0))$ which equals
$\Lambda(\dot\gamma_{\xi_v}(T-T_0))$.

By the fact that $M_1$ and $M_2$ have the same lens data, we
conclude that \[\ell(\phi_{0*}(\dot\gamma_{\xi_v}(T-T_0)))\geq
T_0,\] so the point $\exp(T\phi_{0*}\xi_v)$ does indeed exist and
lie in $M_2$.

Now let $T'$ be another time in between $T_0$ and $T_1$, and let
$\xi_v'$ be the corresponding vector. Without loss of generality we
may assume that $\Delta T= T-T'>0$. Then we have the following
identity

\[
\exp^{g_1}(t\xi_v)=\exp^{g_1}((t-\Delta T)\xi_v').
\]

By the definitions of $T_0$ and $T_1$, the curve
$\exp^{g_1}(t\xi_v):0\leq t\leq \Delta T$ is a geodesic segment
lying entirely within the interior of $U_1$.  Since $\phi_0$ is an
isometry on $U_1$, the vectors $\phi_{0*}\xi_v$ and
$\phi_{0*}\xi_v'$ must also be tangent to a common geodesic at a
distance of $\Delta T$.  Hence

\[
\exp^{g_2}(t\phi_{0*}\xi_v)=\exp^{g_2}((t-\Delta T)\phi_{0*}\xi_v').
\]
Setting $t=T$, we obtain the needed result.
\end{proof}

\begin{proposition}
For fixed $x_0$, $\ti{\phi}(x_0,v)$ is constant.
\end{proposition}
\begin{proof}
The strategy here is to prove that $\ti{\phi}(x_0,v)$ is locally
constant.  Then the statement follows from the connectedness of the
sphere.  First, we need a lemma.

For a pair of points in $\ti{M}$, let $d(x,y)$ denote the distance
between them.  This function is not necessarily smooth, even off the
diagonal.  However, the next lemma shows that the squared distance
function $d(x,y)^2$ is as smooth as the metric for $(x,y)$
sufficiently close to each other.

\begin{lemma}
Let $\ti{M}$ be as above (with subscript omitted).  For every $x_0$
in the interior of $\ti{M}$, there exists a positive number $r$ such
that the squared distance function is analytic on the set
\[\Delta_r(x_0)=\{(x,y):d(x,x_0)<r, d(x,y)<r\}.\]
If $K$ is a compact set contained within the interior of $\ti{M}$,
then there is an open $\sc{O}\subset \ti{M}$ containing $K$ and a
positive number $r$ such that the squared distance function is
analytic on the set
\[\Delta_{\sc{O},r}(K)=\{(x,y):x\in\sc{O}, d(x,y)<r\}.\]
\end{lemma}
\begin{proof}
We choose $r>0$ so that the ball $B_{2r}(x_0)$ is contained within
$\ti{M}$ and is geodesically convex (see \cite{St}, Theorem 6.2,
noting that the restriction on the radius is only that it be
sufficiently small). By definition every pair of points in
$B_{2r}(x_0)$ is joined by a unique geodesic segment contained
entirely within $B_{2r}(x_0)$. Moreover, that segment is
length-minimizing.

Now assume that $(x,y)\in \Delta_r(x_0)$.  Then there is exactly one
geodesic segment connecting them whose length is less than $r$.
Indeed, there is at least one, since the two points lie within
$B_{2r}(x_0)$.  If there were another geodesic segment connecting
them, it would have to leave $B_{2r}(x_0)$ at some point.  Since
$d(x_0,x)<r$, such a segment would necessarily have length greater
than $r$.

This shows that the mapping
\[\{(x,v):d(x,x_0)<r,|v|_g<r\}\rightarrow \Delta_r(x_0)\]
given by $(x,v)\mapsto (x,\exp_x(v))$ is bijective.  Naturally, the
exponential map is analytic as long as the metric is analytic.  By
the inverse function theorem, $\exp^g$ gives a diffeomorphism
between these two sets. Through this diffeomorphism, the squared
distance function is expressed $d(x,y)^2=g_{ij}v^iv^j$, which is
analytic as long as $g$ is analytic.

The second statement of the lemma follows from the first by a
compactness argument.  Indeed, for every $x_0\in K$ we take the
number $r$ from the first statement and form the ball $B_r(x_0)$.
All such balls form an open cover of $K$. We take a finite subcover,
let $\sc{O}$ be the union of its members, and let $r$ be the
smallest radius in that subcover.
\end{proof}

Now fix a vector $(x_0,v)$ and choose $T=T(v)$ and $\xi_v$.  Let
$y_0\in U_1$ be the basepoint of the vector $\xi_v$ so that
$x_0=\exp^{g_1}_{y_0}(T\xi_v)$. Also, let $\gamma_1=\gamma_{x_0y_0}$
denote the geodesic segment connecting the two points.  Let
$\sc{O_1}$, $r_1$ be the open set and radius corresponding to the
compact set $\gamma_1$ as in the lemma.

In $\ti{M_2}$, we have a corresponding segment $\gamma_2$ between
the points $\phi_0(y_0)$ and $\ti{\phi}(x_0,v)$.  It is given by the
curve
\[\exp^{g_2}(t\phi_{0*}\xi_v):,\,\,0\leq t \leq T.\]
Let $r_2$ be the radius corresponding to $\gamma_2$ as in the lemma.

 Let $r$ denote the positive number
 \[ r=\inf\{d(\gamma_1,\p \ti{M_1}), d(y_0,\p U_1), r_1,r_2\}.\]

By continuity, there exists a neighborhood $N$ of $v$ in
$S_{x_0}M_1$ sufficiently small such that for all $w\in N$,
\[d(\gamma_{x_0,-v}(t),\gamma_{x_0,-w}(t))<r\] for all $t$ in the
interval $[0,T]$.  The restrictions on $r$ guarantee that the curve
$\gamma_{x_0,-w}(t)$ remain within $\ti{M_1}$ and that its endpoint,
$\gamma_{x_0,-w}(T)$, be in the interior of $U_1$. For each $w$, let
$\eta_w=-\dot\gamma_{x_0,-w}(T)$.  We prove
$\ti{\phi}(x_0,v)=\ti{\phi}(x_0,w)$ by breaking this into the two
equations:
\begin{equation}\label{puppydog}
\exp^{g_2}(T\phi_{0*}\xi_v)=\exp^{g_2}(T\phi_{0*}\eta_w);
\end{equation}
\begin{equation}\label{hounddog}
\exp^{g_2}(T\phi_{0*}\eta_w)=\exp^{g_2}(T(w)\phi_{0*}\xi_w).
\end{equation}
See Figure \ref{locally constant}.

\begin{figure}[htp]
\centering
\includegraphics{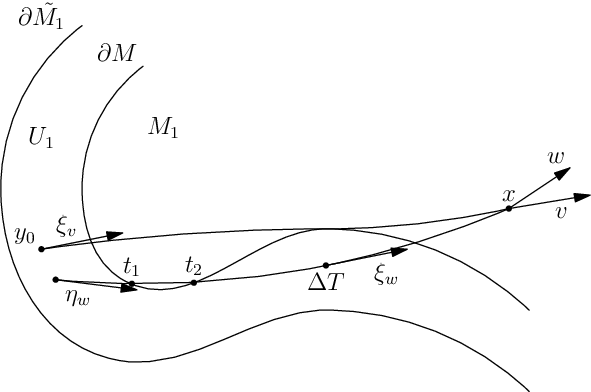}
\caption{} \label{locally constant}
\end{figure}

Consider the function
\[\rho_1(t)=d^2_{g_1}(\exp(t\xi_v),\exp(t\eta_w)).\]
By our choices of $r$ and $N$, and by the lemma, this is an analytic
function for $t\in [0,T]$.

Now we consider $\ti{M_2}$, and define
\[\rho_2(t)=d_{g_2}^2(\exp^{g_2}(t\phi_{0*}\xi_v),\exp^{g_2}(t\phi_{0*}\eta_w)).\]
Since $\phi_0|_{U_1}:U_1\ra U_2$ is an isometry, the functions
$\rho_1$ and $\rho_2$ must coincide for small values of $t$. Also,
we note that the function $\rho_2(t)$ is analytic as long as
$d_{g_2}(\exp^{g_2}(t\phi_{0*}\xi_v),\exp^{g_2}(t\phi_{0*}\eta_w))<r$,
since $r$ was chosen to be smaller than $r_2$. Therefore, by
analytic continuation, the functions $\rho_1$ and $\rho_2$ are equal
up to the first point $t_0$ where $\rho_2=r^2$. But by continuity,
we would then have $\rho_1(t_0)=r^2$, which does not occur.
Therefore, we see that $\rho_1(t)=\rho_2(t)$ throughout the interval
$0\leq t \leq T$.  In particular, we find that $\rho_2(T)=0$, which
verifies equation \eqref{puppydog}.

If $T$ lies between $T_0(w)$ and $T_1(w)$, then equation
\eqref{hounddog} is nothing but a restatement that
$\ti{\phi}(x_0,w)$ is well defined. Clearly, $T>T_0(w)$, so assume
that it is also greater than $T_1(w)$ (The possiblity that
$T=T_1(w)$ is ruled out by the fact that $\gamma_{x,-w}(T_1(w))\in
\p M$). Fix a number $T'=T(w)$, and a corresponding $\xi_w$. If we
let $\Delta T=T-T'$, then we have the equation
\[\exp^{g_1}(t\eta_w)=\exp^{g_1}((t-\Delta T)\xi_w).\]

Let $\gamma(t)=\exp^{g_1}(t\eta_w)$, for $0\leq t\leq \Delta T$. If
it happens that it lies entirely within $U_1$, then the same proof
that we used to show that $\ti{\phi}$ is well-defined will verify
equation \eqref{hounddog}.  So assume that $\gamma(t)$ does not lie
entirely within $U_1$. Then we can uniquely partition it into
subsegments which alternately lie in $U_1$ and $M_1$. Indeed, we
define
\[\begin{array}{cc} t_0= & 0, \\
t_1= & \inf \{t>0:\gamma(t)\notin U_1\}, \\
t_2= & \inf \{t> t_1:\gamma(t)\notin M_1\},\\
\vdots & \vdots \\
t_{m-1}= & \inf \{t> t_{m-2}:\gamma(t)\notin M_1\}, \\
t_m= & \Delta T.\end{array}\] The fact that
$0=t_0<t_1<...<t_m=\Delta T$ is true follows from the same reasoning
that was used above to prove that $T_1(v)>T_0(v)$.  The fact that
this partition is finite also follows from the analytic assumption.
Indeed an analytic curve segment cannot intersect the analytic $\bo$
more than a finite number of times without being entirely contained
within it. Note that the segment $\gamma|_{[t_k,t_{k+1}]}$ lies in
$U_1$ or $M_1$ according to whether $k$ is even or odd,
respectively.  In particular, $m$ is odd.

To prove equation \eqref{hounddog}, we will use induction to show
that for all $k=1,2,\dots,m$, and all $t\in [0,T]$,
\begin{equation}\label{wienerdog}
\exp^{g_2}(t\phi_{0*}\dot\gamma(0))=\exp^{g_2}((t-t_k)\phi_{0*}\dot\gamma(t_k)).
\end{equation}
Then setting $k=m$ and $t=T$ yields the result.

Step 1: $\eta_w=\dot\gamma(0)$ and $\dot\gamma(t_1)$ lie on the
geodesic $\gamma$ at a distance of $t_1$ from each other.  Since
this segment lies within $U_1$ and since $\phi_0$ is an isometry of
$U_1$ to $U_2$, we see that $\phi_{0*}\dot\gamma(0)$ and
$\phi_{0*}\dot\gamma(t_1)$ also lie on a common geodesic at the same
distance apart.  Hence equation \eqref{wienerdog} is established for
$k=1$.

Step 2: The next segment of $\gamma$ lies within $M_1$.  Indeed we
have the following:
\[\dot\gamma(t_2)=\Sigma_{g_1}(\dot\gamma(t_1)),\,\,\ell_{g_1}(\dot\gamma(t_1))=t_2-t_1.\]
Since $M_2$ has the same lens data as $M_1$, we see that
$\phi_{0*}\dot\gamma(t_1)$ and $\phi_{0*}\dot\gamma(t_2)$ are
connected by a geodesic across $M_2$ with the same length $t_2-t_1$.
Together with step 1, this shows that $\phi_{0*}\dot\gamma(t_2)$
lies tangent to the same geodesic as $\phi_{0*}\dot\gamma(0)$ at a
distance of $t_2$.  Hence Equation 8 is established for $k=2$.

Step 3: By induction, we may repeat these steps, establishing
equation \eqref{wienerdog} for all $k$ up to $k=m$.

\end{proof}

For $x\in M_1$, set $\phi(x)=\ti{\phi}(x,v)$.  If $x\in V_1$, then
we are in the domain of the boundary normal coordinates.  We choose
$v=\frac{\p}{\p x^n}$.  Then $\gamma_{x,-v}$ is the geodesic segment
from $x$ to $U_1$ normal to $\bo$.  We choose $T(v),\xi_v$ so that
$x=\exp^{g_1}(T(v)\xi_v)$ and note that the segment
$\exp^{g_1}(t\xi_v): 0\leq t\leq T_v$ lies entirely in $V$.  We have
\[\phi_0(x)=\exp(T(v)\phi_{0*}\xi_v)=\ti{\phi}(x,v)=\phi(x)\]
The first equation is true by the fact that $\phi_0$ is an isometry
on $V_1$.  Hence $\phi$ and $\phi_0$ agree on their common domains.
Gluing them together, we form
\[\phi:\ti{M_1}\ra \ti{M_2}.\]

\section{That $\phi$ is a diffeomorphism and an isometry}

Let $(x,v)$ lie in the interior of $S_xM_1$.  We choose $T=T(v)$ and
$\xi_v\in SU_1$ so that $\gamma(t)=\exp^{g_1}(t\xi_v)$ is a geodesic
in $\ti{M_1}$ that reaches $x$ at time $T$.  In the segment
$[0,T-T_0]$, $\gamma$ lies entirely within $U_1$; whereas on the
segment $[T-T_0,T+\delta]$, it lies entirely within $M_1$ for
$\delta$ sufficiently small. For the first segment, we see that
\[\phi(\gamma(t))=\phi_0(\gamma(t))=\exp^{g_2}(t\phi_{0*}\xi_v);\]
where the second equality holds by the fact that $\phi_0$ is an
isometry on $U_1$.

On the second segment, we see that for each pair
$(\gamma(t),\dot\gamma(t))$, we can choose the same $\xi_v$ for
$\xi_{\dot\gamma(t)}$ with the corresponding $T(\dot\gamma(t))=t$.
So, for all $t\in [0,T+\delta]$,

\begin{equation}\label{bulldog}
\phi(\gamma(t))=\exp^{g_2}(t\phi_{0*}\xi_v).
\end{equation}
This can be rewritten in the form:
\begin{equation}\label{maddog}
\phi(\gamma(t))=\exp^{g_2}(t\phi_{0*}\dot\gamma(0)).
\end{equation}

In fact, the latter equation is true for any geodesic segment
$\gamma|{[0,T]}$ that can be partitioned into $\gamma|{[0,a]}$ and
$\gamma|{[a,T]}$ with the two subsegments lying entirely in $U_1$
and $M_1$ respectively.

\begin{proposition}
$\phi:\ti{M_1}\ra \ti{M_2}$ is bijective.
\end{proposition}
\begin{proof}
Reversing the roles of $\ti{M_1}$ and $\ti{M_2}$, we can define a
mapping $\psi:\ti{M_2}\ra\ti{M_1}$ by the same process by which we
defined $\phi$.  In particular it would extend $\phi_0^{-1}$.

The analog to equation \eqref{maddog} is:
\begin{equation}\label{hangdog}
\psi(\beta(t))=\exp^{g_1}(t\phi^{-1}_{0*}\dot\beta(0)),
\end{equation}
where $\beta$ is any geodesic segment composed of two subsegments
contained in $U_2$ and $M_2$, respectively.

Using the notation from above, we would like to prove that
$\psi\circ\phi(\gamma(t))=\gamma(t)$ for $t\in[0,T]$.  To that end,
we will first show that the geodesic segment
$\beta(t)=\phi(\gamma(t))$ is of the type valid for equation
\eqref{hangdog}.

We note that for $t\in[0,T-T_0]$, $\gamma(t)\in U_1$ so
$\phi(\gamma(t))$ must lie in $U_2$.  For $t\in[T-T_0,T]$,
$\gamma(t)\in M_1$, so by Proposition 1, $\phi(\gamma(t))\in M_2$.

Therefore, we may apply equation \eqref{hangdog}, which yields:
\[\begin{array}{ccc}
\psi\circ\phi(\gamma(t)) & = &
\exp^{g_1}(t\phi_{0*}^{-1}\dot\beta(0))
\\
 & = & \exp^{g_1}(t\phi_{0*}^{-1}\phi_{0*}\dot\gamma(0)) \\
 & = & \exp^{g_1}(t\dot\gamma(0)) \\ & = & \gamma(t).\end{array}\]

Since every point in $M_1$ lies on some such curve $\gamma(t)$, we
conclude that $\psi\circ\phi=id$ on $M_1$. But we know that the same
identity is true on $V_1$, so it is true on all of $\ti{M_1}$. By
the symmetry of the construction, we conclude $\phi\circ\psi$ is
also the identity.
\end{proof}

\begin{proposition}
$\phi:\ti{M_1}\ra \ti{M_2}$ is an analytic isometry.
\end{proposition}
\begin{proof}
Since $\phi$ is bijective, it is sufficient to prove the statement
locally. These properties are already known on $V_1$ where
$\phi=\phi_0$, so we assume $x$ lies in the interior $M_1$. First we
show that all directional derivatives of $\phi$ exist.  Indeed, by
differentiating equation \eqref{bulldog} at $t=T(v)$, we obtain:
\[D_v\phi(x)=\p_t(\phi\circ\gamma)(T(v))=\p_t
\exp^{g_2}(t\phi_{0*}\xi_v)|_{t=T(v)}.
\]
Clearly, the quantity on the right side exists.  What's more, it is
a vector of length $1$.  We conclude that $\phi_*$ exists and
preserves lengths of vectors.  In particular it is nonsingular.

From the fact that it preserves lengths, we derive smoothness.
Indeed, $g_1=\phi^*g_2$, which has the coordinate form:
\[g_1(x)_{ij}=\phi^k_{,\,i}(x)\,g_2(\phi(x))_{kj}.\]
This yields:
\[g_1(x)_{ij}g_2(\phi(x))^{jl}=\phi^k_{,\,i}(x).\]
The left side is once differentiable; hence $\phi$ is twice
differentiable.  But then that implies the left side to be twice
differentiable which shows $\phi$ to be three-times differentiable.
By an obvious application of induction, $\phi$ must be smooth.

It only remains to prove that $\phi$ is analytic.  Of course this is
already known in $V_1$. Since $\phi$ is a smooth isometry, we can
state
\[\phi(\exp^{g_1} \xi)=\exp^{g_2}(\phi_*\xi)\]for any vector $\xi\in
T\ti{M_1}$. Given $x_0$ in the interior of $M_1$, consider a ball of
radius $r$, centered at $x_0$, which is strictly geodesically
convex, and choose any point $y_0$ within this ball. Then there is a
unique $\xi_0\in T_yM_1$ satisfying $|\xi_0|_{g_1}<r$ and
$\exp_{y_0}\xi_0=x$. Moreover, $x_0$ and $y_0$ are not conjugate
along the corresponding geodesic, so
\[\exp_{y_0}^{g_1}:\xi\mapsto x\] is a local diffeomorphism about
$\xi_0$.  It is analytic by the fact that $g_1$ is analytic.
Consequently, it's inverse function is analytic. Let $\xi(x)$ denote
the inverse, which is defined for $x$ in some neighborhood of $x_0$.
Then we see that
\[\phi(x)=\exp^{g_2}_{\phi(y_0)}(\phi_*\xi(x)).\]  $y_0$ is constant, so
$\phi_*$ is a constant linear mapping.  From the fact that $g_2$ is
analytic, we see that this mapping is also analytic.
\end{proof}


\begin{thebibliography}
\frenchspacing

\bibitem[BI]{BI} {\sc D. Burago and S. Ivanov}, Boundary rigidity and filling
volume minimality of metrics close to a flat one, manuscript, 2005.

\bibitem[BCG]{BCG}{\sc G. Besson, G. Courtois, and S. Gallot},
Entropies et rigidit\'es des espaces localement sym\'etriques de
courbure strictment n\'egative, {\it Geom. Funct. Anal.}, {\bf
5}(1995), 731-799.

\bibitem[C1]{C1} {\sc C. Croke}, Rigidity and the distance between boundary
points, {\it J. Differential Geom.}, {\bf 33}(1991), no. 2,
445--464.

\bibitem[C2]{C2} {\sc C. Croke}, Boundary and lens rigidity of finite quotients, {\it Proc. AMS}, {\bf 133}(2005), no. 12, 3663--3668.

\bibitem[CDS]{CDS} {\sc C. Croke, N. Dairbekov, V. Sharafutdinov}, Local
 boundary rigidity of a compact Riemannian manifold with curvature
 bounded above, {\it Trans. Amer. Math. Soc.} {\bf 352}(2000), no. 9, 3937--3956.

\bibitem[CK]{CK} {\sc C. Croke and B. Kleiner}, Conjugacy and rigidity for manifolds with a parallel vector field, {\it J. Diff. Geom.} 39(1994), 659--680.

\bibitem[E]{E} {\sc G. Eskin}, Inverse scattering problem in anisotropic media, {\em  Comm. Math. Phys.} {\bf 199}(1998), no. 2, 471--491.

\bibitem[Gr]{Gr} {\sc M. Gromov}, Filling Riemannian manifolds, {\it J. Diff. Geometry} {\bf 18}(1983), no. 1, 1--148.

\bibitem[LSU]{LSU} {\sc M. Lassas, V. Sharafutdinov and G. Uhlmann}, Semiglobal boundary rigidity for Riemannian metrics, {\it Math. Ann.}, {\bf 325}(2003), 767--793.

\bibitem[Mi]{M} {\sc R. Michel}, Sur la rigidit\'e impos\'ee par la
longueur des g\'eod\'esiques, {\it Invent. Math.} {\bf 65}(1981),
71--83.

\bibitem[PoR]{PoR} {\sc M. Porrati and R. Rabadan}, Boundary rigidity and holography,
{\it J. High Energy Phys.} {\bf 0401}, 034-057.

\bibitem[Po]{Po}
{\sc M.M. Postnikov}, Geometry VI: Riemannian Geometry (2001).

\bibitem[PU]{PU} {\sc L. Pestov and G. Uhlmann}, Two dimensional simple
compact manifolds with boundary are boundary rigid,   {\it Ann.
Math.} {\bf 161}(2)(2005), 1089--1106.

\bibitem[SU1]{SU1}
{\sc P. Stefanov and G.Uhlmann}, Rigidity for metrics with the same
lengths of geodesics, {\em Math. Res. Lett.} {\bf 5}(1998), 83--96.

\bibitem[SU2]{SU2}
{\sc P. Stefanov and G.Uhlmann}, Boundary rigidity and stability for
generic simple metrics,  {\em J. Amer. Math. Soc.} {\bf 18}(2005),
975--1003.

\bibitem[SU3]{SU3}
{\sc P. Stefanov and G.Uhlmann}, Local Lens Rigidity with Incomplete
Data for a class of non-Simple Riemannian Manifolds.

\bibitem[St]{St}
{\sc S. Sternberg}, Lectures on Differential Geometry (1964).


Cauchy data,  {\it IMA publications} {\bf Vol 137}, ``Geometric
methods in inverse problems and PDE control" (2003), 263-288.


\end{thebibliography}
\end{document}